\documentclass{aptpub}

\authornames{BIAO WU}

\shorttitle{Multiagent models in time-varying and random
environment}

\usepackage{amstext, graphics, amsmath, latexsym, mathrsfs, amssymb, amscd, amsfonts}

\newcommand{\bin}[2]{\mbox{$\left({}^{#1}_{#2}\right)$}}

\usepackage{natbib}

\numberwithin{equation}{section}

\begin{document}

\title{MULTIAGENT MODELS IN TIME-VARYING \\ AND RANDOM ENVIRONMENT}

\authorone[Carleton University]{BIAO WU}

\addressone{School of Mathematics and Statistics, Carleton
University, 1125 Colonel By Drive, Ottawa, ON, Canada K1S 5B6 Email
address: biaow@math.carleton.ca}

\begin{abstract}
In this paper we study multiagent models with time-varying type
change. Assume that there exist a closed system of $N$ agents
classified into $r$ types according to their states of an internal
system; each agent changes its type by an internal dynamics of the
internal states or by the relative frequency of different internal
states among the others, e.g., multinomial sampling. We investigate
the asymptotic behavior of the empirical distributions of the
agents' types as $N$ goes to infinity, by the weak convergence
criteria for time-inhomogeneous Markov processes and the theory of
Volterra integral equations of the second kind. We also prove
convergence theorems of these models evolving in random environment.

\end{abstract}



\section{INTRODUCTION}\label{Introduction}

Agent-based models (ABMs), or multiagent systems (MAS), arise from
many areas of science and social sciences such as ecology,
artificial intelligence, communication networks, sociology,
economics; see e.g. Ferber (1999, 1995), Ouelhadj (1996), and
Wooldridge (1995). Goldstone and Janssen (2006) studied the
collective behavior of the agent-based computational models, which
build social structures from the `bottom-up'. We give some of
 the attractive features of ABMs presented in their paper.
 First, ABMs can describe precise mathematical formulation, which
make clear, quantitative and objective predictions possible. Second,
ABMs can bridge the explanations that link the analysis of the
individual agent level and the analysis of the emergent group level.

In this paper we will focus on the agent based modeling in social or
economic discipline. F\"{o}llmer and Schweizer (1993) considered an
interacting agent financial model in which they used Black's (1986)
classification of traders: information traders and noise traders.
Lux (1995, 1997, 1998) studied a model of three types of traders
which can probabilistically change their types. Horst (2000, 2001,
2002, 2005) kept some aspects of F\"{o}llmer and Schweizer's model
and considered interacting agent models with local and global
interactions. Horst assumed that the set of active agents
$\mathcal{A}$ is countable and there is a sequence of finite sets
$\mathcal{A}_n$ satisfying $\lim_{n\to
\infty}\mathcal{A}_n=\mathcal{A}$. In Horst's example, the traders
are divided into fundamental traders and noise traders, and the
fundamental traders are divided into optimistic and pessimistic
fundamental traders. Horst introduced the concept of individual mood
into his models. At each period $t\ge 0$, each fundamental trader
has its own mood, e.g., $x^a_t=+1$ or $x^a_t=-1$, that is to say,
the fundamental trader is an optimist or a pessimist. Let $C$ be a
fixed set of individual states, i.e., $x^a_t\in C$ for each $a\in
\mathcal{A}$ and $t\ge 0$. Let $x_t=\{x^a_t\}_{a\in \mathcal{A}}$.
Horst defined the empirical distribution, which is called {\em mood
of the market}, as follows:
$$\rho_t=\rho(x_t)=\lim_{n\to \infty}\frac{1}{|\mathcal{A}_n|}\sum_{a\in
\mathcal{A}_n}\delta_{x^a_t}(\cdot).$$ The market mood is one of the
main driving forces of Horst's interacting agent models.

So far, we have seen the importance of the empirical distribution of
the individual states which links the behavior of individual agent
level and the emergent laws of collective level. The multiagent
models of this paper arise from social or economic background; some
features of them show similarities with the Wright-Fisher model in
population biology , and the Voter model, see e.g. Either and Kurtz
(1986), and Holley (1975). Instead of giving the precise definition
of the agents, we describe the properties and behaviors of the
agents rather intuitively. The multiagent models here share some
similarities with ABMs in other disciplines. The most general
assumptions about the mechanism of the multiagent models are as
follows:
\begin{enumerate}
  \item The time is in nonnegative integer units, denoted by $k\ge 0$.

  \item There are fixed $N\ge 2$ agents in the multiagent system at all times. There are no entries of new agents into the system or
exits of current agents from the system.

  \item There is an internal system with which all agents are
  concerned. The internal system has $r\ge 2$ states which we simply
  denote by $1$, $\cdots$, $r$. The internal system will not change with time $k$.
  That is to say, at any time $k\ge 0$, there is no new state added to the
  internal system and there is no existing state removed from it.
  The behavior of internal states are observed by all agents. Each agent has one and only one internal state at each time $k\ge
  0$. Thus, the agents are classified into $r$ types, according to their
  internal states.

  \item Assume that $n^N_{i}(k)$ ($1\le i \le r$) is the number of agents of type $i$ at time
$k$ and $\mathbf{n}^N(k)=(n^N_{1}(k), n^N_{2}(k),\cdots,
n^N_{r}(k))$ is the distribution of all agents among the $r$ types.
By the second assumption, $n^N_1(k)+ \cdots+ n^N_r(k)=N$ and
$\frac{\mathbf{n}^N(k)}{N}$ is the empirical distribution of the
types at $k\ge 0$.

\item Assume that $\{(p_{N, i, j}(k))_{r\times r}, k\ge 0\}$ is a sequence of
deterministic {\em stochastic-matrix} valued functions which
represent the external environment of the multiagent system.

  \item Based on all the information of the agents' types and the environment up to time $k$, each agent has an independent strategy of
   probabilistically choosing its type for the next time unit $k+1$. The strategy of an agent is realized by keeping or changing
  its type. The agents of the same type have a common strategy. That is to say, from time $k$ to $k+1$, the
  agents of type $i$ switch to type $j$ with probability $p_{N, i,
  j}(k)$. This process of changing types
  occurs locally among agents of the same type, and it is called {\em internal dynamics}.

  \item From time $k$ to $k+1$, there also occurs another process of global type
 change. When we make this assumption, we would
  make a minor change on the fifth assumption, to say
  that the internal dynamics occurs from time $k$ to $k+\frac{1}{2}$, rather than from time $k$ to $k+1$.
   Based on all the information up to time $k+\frac{1}{2}$, each agent independently determines its type by
$\frac{\mathbf{n}^N(k+\frac{1}{2})}{N}$. That is to say, for any
agent, regardless of its type at time $k+\frac{1}{2}$, the
probability of its new type being $i$ at time $k+1$ is
$\frac{n^N_i(k+\frac{1}{2})}{N}$. Therefore,
$
{\mathbf n}^N(k+1) \sim \mbox{multinomial}\biggl(N,
\frac{\mathbf{n}^N(k+\frac{1}{2})}{N}\biggl).
$

\end{enumerate}

Note that we can change the number $\frac{1}{2}$ in the seventh
assumption by any number $c$ ($0<c<1$). The external environment in
the fifth assumption can be external economic fundamentals. An
example of this is given by Example \ref{examp:EconFund}. Based on
the assumptions 1-6, we can construct the multiagent model with
internal dynamics (MAMWID). If we have all the seven assumptions, we
can construct multiagent model with internal dynamics and
multinomial sampling (MAMWIDAMS). The multinomial sampling is a kind
of interaction among the agents. When we assume that $\{(p_{N, i,
j}(k))_{r\times r}, k\ge 0\}$ is a random sequence, we can construct
multiagent model with internal dynamics and random environment
(MAMWIDARE); and multiagent model with internal dynamics,
multinomial sampling, and random environment (MAMWIDAMSARE)
respectively.

In this paper, we will mainly study the asymptotic behaviors of the
empirical distribution, $\{\frac{\mathbf{n}^N([Nt])}{N}, t\ge 0\}$,
of the types, as $N\to \infty$. Lux assumed the number of the agents
to be finite, but he didn't consider the asymptotics of the
empirical distribution of the types as $N$ becomes large.
F\"{o}llmer, Schweizer, and Horst assumed that the number of the
agents is countable; they do not have the question we discuss here.
Another feature of our multiagent models is that the transition
structure of the internal dynamics is time-inhomogeneous.

This paper is the first attempt of a systematic study of the
interacting agent financial systems. Another working paper of the
author, which goes one step further than this one, focuses on the
interacting agent feedback finance models, see Wu (2006).

This paper is organized as follows. In Subsection
\ref{sse:MAMIDTVE}, we formulate MAMWID and MAMWIDAMS; and state
their convergence in Theorem \ref{thm:MAMIDTVE}. In Subsection
\ref{sse:MAMSIRE}, we formulate MAMWIDARE and MAMWIDAMSARE; and
state their joint and annealed convergence in Theorem
\ref{thm:JointAnnealedMarginalRanCon}. In Section
\ref{se:ProofOfMAMWID}, we prove Theorem \ref{thm:MAMIDTVE}; and in
Section 4, we prove Theorem \ref{thm:JointAnnealedMarginalRanCon}.
In Appendix \ref{se:WCCFTInhomMP}, we state weak convergence
criteria for time inhomogeneous Markov processes, which are the main
tools of this paper.

\section{Formulation of the Multiagent Models and Main
Results}\label{se:FormResults}

\subsection{Multiagent models in deterministic time-varying environment}\label{sse:MAMIDTVE}

Now we formulate the multiagent model with internal dynamics
(MAMWID) based on the assumptions 1-6 in Section \ref{Introduction}.
Let $ A(t)=(a_{i,j}(t))_{r\times r}$ be a $r\times r$
 matrix-valued $c\grave{a}dl\grave{a}g$ function on $[0, \infty)$, which satisfies the conditions
\begin{enumerate}
 \item[1)] For each $t\ge 0$, $A(t) \mathbf{e}={\mathbf 0}_{r\times 1}$, where $\mathbf{e}=[1, \cdots, 1]'$ and ``$'$'' denotes transpose;
 \item[2)] For each $t\ge 0$, $1\le i,j\le r$, $i\not=j$, $a_{i,j}(t)\ge 0$.
\end{enumerate}
Fix $N\ge 1$, let $ A_N(t)=(a_{N,i,j}(t))_{r\times r}$ be a $r\times
r$
 matrix-valued function on $[0, \infty)$, which satisfies the conditions
\begin{enumerate}
 \item[1)] For each $t\ge 0$, $A_N(t) \mathbf{e}={\mathbf 0}_{r\times 1}$;
 \item[2)] For each $t\ge 0$, $1\le i,j\le r$, $i\not=j$, $a_{N,i,j}(t)\ge
 0$;
 \item[3)] $A_N$ is a stepwise function, i.e., $A_N(t)$ is a constant on $[\frac{k}{N},
\frac{k+1}{N})$ for each $k\ge 0$.
\end{enumerate}
Therefore $A(t)$ and $\{A_N(t)\}$ are {\em Q-matrix} valued
functions.

Let $R^{r\times r}$ be the Euclidean space corresponding to the
$r\times r$ square matrix. For each $N\ge 1$, and $k\ge 0$, let
\begin{equation}\label{eq:PNANk}
P_{N,k}=(p_{N, i, j}(k))_{r\times r}=I+\frac{1}{N}A_N(\frac{k}{N}),
\end{equation}
where $I$ is the identity matrix of order $r$, and $p_{N, i, j}(k)$
is the probability of each agent of type $i$ switching to type $j$
at time $k+1$. The definition of $P_{N,k}$ is valid since for large
enough $N$, $P_{N,k}$ is a stochastic matrix, which we call internal
transition matrix of MAMWID.

We are ready to formulate MAMWID. For $k\ge 0$, the transition
between ${\mathbf n}^N(k)$ and ${\mathbf n}^N(k+1)$ is determined by
the sixth assumption in Section \ref{Introduction} as follows. For
$1\le i\le r$, each agent of type $i$ can change its type to $j$,
with probability $p_{N, i, j}(k)$ ($1\le j\le r$). Since the
$n^N_{i}(k)$ agents independently make their transitions, the
distribution of these $n^N_{i}(k)$ agents among the $r$ types at
time $k+1$ is a random vector denoted by $\mathbf{\Xi}_{N, k,
i}=(\xi_{N, k, i, 1}, \cdots, \xi_{N, k, i, r})$, which satisfies
\begin{equation}\label{eq:XiNki}
\mathbf{\Xi}_{N, k, i} \sim \mbox{multinomial}(n^N_{i}(k),
P_{N,k,i\cdot}),
\end{equation} where $P_{N,k,i\cdot}$ is the $i$-th
row of the matrix $P_{N,k}$. Since agents in different type change
their types independently, $\mathbf{\Xi}_{N,k,1}, \cdots,
\mathbf{\Xi}_{N,k,r}$ are independent. The distribution of all the
agents at time $k+1$ is
\begin{equation}\label{eq:nk+In}
{\mathbf
n}^N(k+1)\equiv\mathbf{\Xi}_{N,k,1}+\cdots+\mathbf{\Xi}_{N,k, r}.
\end{equation}
 The sequence $\{ {\mathbf n}^N(k), k \ge 0\}$ defined this way
is a time inhomogeneous Markov chain. At last, we
 define
 \begin{equation}\label{eq:XNt}
 \mathbf{X}^N(t)\equiv\frac{\mathbf{n}^N([Nt])}{N}.
 \end{equation}

We introduce some notations. We put $Z_+=\{0, 1, \cdots\}$ and
$R_+=[0,\infty)$,
\begin{equation*}
K_N=\{N^{-1}\mathbf{\alpha}: \mathbf{\alpha}\in (Z_+)^r, \mbox{
}\sum_{i=1}^{r}\alpha_i=N\},
\end{equation*}
and
\begin{equation*}
K=\{\mathbf{\alpha}: \mathbf{\alpha}\in (R_+)^r, \mbox{
}\sum_{i=1}^{r}\alpha_i=1\}.
\end{equation*}

We define the time-dependent generator
    $\{G_A(t), 0\le t<\infty\} $ on $C^1(K)$: for each $f\in C^1(K)$ and $t\ge 0$,
\begin{equation}\label{eq:genGeneralInhom}
G_A(t)f({\mathbf p})={\mathbf p}A(t)\frac{\partial f}{\partial
\mathbf{x}}'.
\end{equation}
It is clear that, for each $t\ge 0$, $\mathscr{D}(G_A(t))= C^1(K)$
and $\mathscr{D}(G_A)= C^1(K)$ is the common domain of
 the generator $\{G_A(t), 0\le t<\infty\} $, and $D= C^2(K)$ is a subalgebra contained in
 $\mathscr{D}(G_A)$.

Next, we illustrate the transition structure of the multiagent model
with internal dynamics and multinomial sampling (MAMWIDAMS) by the
following diagram:
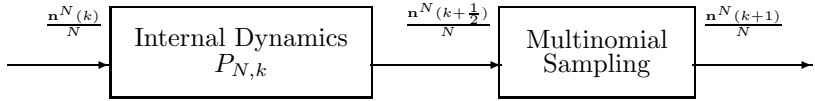
\begin{figure}[h]
\begin{center}
\vspace{1.5cm}
\begin{picture}(250,0)
$\stackrel{\frac{\mathbf{n}^N(k)}{N}}{\put(-25.0,12){\vector(1,0){38.5}}}$\put(3,0){\framebox(98.5,34)[c]{$\stackrel{\mbox{Internal
Dynamics}}{P_{N,k}}$}
}$\stackrel{\hspace{4.0cm}\frac{\mathbf{n}^N(k+\frac{1}{2})}{N}}{\put(27.5,12){\vector(1,0){49.5}}}$
\put(3,0){\framebox(73.5,34)[c]{$\stackrel{\mbox{Multinomial}}{\mbox{Sampling}}$}
}$\stackrel{\hspace{2.8cm}\frac{\mathbf{n}^N(k+1)}{N}}{\put(20.5,12){\vector(1,0){45.5}}}$
\end{picture}
\end{center}
\caption{The transition structure of MAMWIDAMS}
\label{MAMWIDAMS:fig}
\end{figure}

Once $\{ \mathbf{n}^N(k), k\ge 0\}$ is defined by the Figure
\ref{MAMWIDAMS:fig}, we define
 \begin{equation}\label{eq:YNt}
 \mathbf{Y}^N(t)\equiv\frac{\mathbf{n}^N([Nt])}{N}.
 \end{equation}

Next, we define the differential operator $G_B$ on $C^2(K)$ by
\begin{equation}\label{eq:DefnOfB}
 G_B=\frac{1}{2}\sum_{i,j=1}^r
b_{ij}({\mathbf p})\frac{\partial^2}{\partial x_i \partial x_j},
\end{equation}
where $b_{ij}({\mathbf p})=p_i(\delta_{ij}-p_{j})$. We also define
the generator $\{G_{AB}(t), 0\le t<\infty\} $ on $C^2(K)$ by
\begin{equation}\label{eq:defnofGABs}
G_{AB}(t)=G_A(t)+G_B, \mbox{ for any $t\ge 0$},
\end{equation}
whose common domain is
$\mathscr{D}(G_{AB})=\mathscr{D}(G_A)\cap\mathscr{D}(G_B)=C^2(K)$.

Define the metric $d_U$ on $D_{R^{r\times r}}[0,
 \infty)$ as follows:
\begin{equation}\label{eq:defndU}
d_U(\mathbf{x},\mathbf{y})=\int_0^{\infty} e^{-u}\sup_{0\le t\le
u}[\|\mathbf{x}(t)-\mathbf{y}(t)\|_r \wedge 1]du, \mbox{
$\mathbf{x}$, $\mathbf{y}\in D_{R^{r\times r}}[0, \infty)$}.
\end{equation}

{\thm\label{thm:MAMIDTVE} Define $P_{N, k}$ by (\ref{eq:PNANk}). Let
$\mu\in \mathscr{P}(K)$. Assume that $\lim_{N\to \infty}d_U(A_N,
A)=0$, where $d_U$ is defined by (\ref{eq:defndU}).
\begin{itemize}
 \item[\textbf{1})] MAMWID. Define $\{ {\mathbf n}^N(k), k \ge 0\}$ by the internal dynamics. Define
$\mathbf{X}^N$ by (\ref{eq:XNt}).  If
$P(\mathbf{X}^N(0))^{-1}\Rightarrow \mu$, then there exists a unique
solution $\mathbf{X}^{\infty}$ of the $D_K[0,\infty)$ martingale
problem for $(G_A, \mu)$ on $C^2(K)$, and $\mathbf{X}^N\Rightarrow
\mathbf{X}^{\infty}$.

\item[\textbf{2})] MAMWIDMS. Define $\{ {\mathbf n}^N(k), k \ge 0\}$ by the Figure
\ref{MAMWIDAMS:fig}, $\mathbf{Y}^N$ by (\ref{eq:YNt}), and
$G_{AB}(t)$ by (\ref{eq:defnofGABs}). If
$P(\mathbf{Y}^N(0))^{-1}\Rightarrow \mu$, then there exists a unique
solution $\mathbf{Y}^{\infty}$ of the $D_K[0,\infty)$ martingale
problem for $(G_{AB}, \mu)$ on $C^3(K)$, and
$\mathbf{Y}^N\Rightarrow \mathbf{Y}^{\infty}$.
\end{itemize}
 }

\subsection{Multiagent models in random environment}\label{sse:MAMSIRE}

In this subsection, we assume that $\{A_N\}$ and $A$ are random
elements which represent an external random environment. We also
assume that $A$ is $C_{R^{r\times r}}[0,\infty)$-valued. Then, we
need the condition $\lim_{N\to\infty}d(A_N, A)=0$ instead of the
condition $\lim_{N\to\infty}d_U(A_N, A)=0$, where $d$ is the
complete separable Skorohod metric on $D_{R^{r\times r}}[0,\infty)$.

Let $(\mathscr{L}, d_U)$ be the subspace of $(D_{R^{r\times
r}}[0,\infty),d_U)$ such that each element of $\mathscr{L}$
satisfies the conditions at the beginning of subsection
\ref{sse:MAMIDTVE}. Let $(\mathscr{L}_c, d_U)$ be the subspace of
$(C_{R^{r\times r}}[0,\infty),d_U)$ such that
$\mathscr{L}_c=\mathscr{L}\bigcap C_{R^{r\times r}}[0,\infty)$. Let
$(\mathscr{P}(D_K[0,\infty)),\rho)$ be the space of probability
measures on $D_K[0,\infty)$ where $\rho$ is the Prohorov metric on
$\mathscr{P}(D_K[0,\infty))$. Then Theorem \ref{thm:MAMIDTVE} shows
that for any $A\in \mathscr{L}_c$, there exist unique $P_A\in
\mathscr{P}(D_K[0,\infty))$ such that under $P_A$, the coordinate
process $Z$ on $D_K[0,\infty)$ is the unique solution of the
martingale problem for the generator $\{G_A(t), 0\le t<\infty\}$ on
$C^2(K)$, and $P_{AB}\in \mathscr{P}(D_K[0,\infty))$ such that under
$P_{AB}$, $Z$ is the unique solution of the martingale problem for
the generator $\{G_{AB}(t), 0\le t<\infty\}$ on $C^3(K)$. We define
$\Phi: \mathscr{L}_c\mapsto \mathscr{P}(D_K[0,\infty))$ by
$\Phi(A)=P_A$ and $\Psi: \mathscr{L}_c\mapsto
\mathscr{P}(D_K[0,\infty))$ by $\Phi(A)=P_{AB}$.

For each $N\ge 1$, let $(\mathscr{L}_N, d_U)$ be the subspace of
$(\mathscr{L}, d_U)$ such that for each $A_N\in \mathscr{L}_N$,
$A_N(t)$ is a constant on $t\in [\frac{k}{N}, \frac{k+1}{N})$ for
each $k\ge 0$. Then for given $A_N\in \mathscr{L}_N$, there are
unique probability measure
     $P_{A_N}\in \mathscr{P}(D_{K_N}[0,\infty))$ which is related to MAMWID, and unique probability measure
      $P_{A_N B}\in \mathscr{P}(D_{K_N}[0,\infty))$ which is related to MAMWIDAMS. Under $P_{A_N}$ or $P_{A_N B}$,
       the coordinate process $Z_N$ on $D_{K_N}[0,\infty)$ has the same distributions as those defined for
        $\mathbf{X}^N$ of MAMWID or $\mathbf{Y}^N$ of MAMWIDAMS. Then we can define $\Phi_N: \mathscr{L}_N\mapsto \mathscr{P}(D_{K_N}[0,\infty))$ and
         $\Psi_N: \mathscr{L}_N\mapsto \mathscr{P}(D_{K_N}[0,\infty))$ correspondingly.

\begin{example}\label{examp:EconFund} Assume that $(H, d_H)$ is a polish space. $\{h^N(\frac{k}{N}), k\ge 0
\}$ is a sequence of external economic fundamentals taking values in
$(H, d_H)$. We assume also that $F$ is a continuous mapping from
$(D_H[0,\infty), d_{S_H})$ onto $(\mathscr{L}, d_U)$, where
$d_{S_H}$ represents the complete separable Skorohod metric on
$D_H[0,\infty)$. Let $A_N(t)=F(h^N(\frac{[Nt]}{N}))$, for any $N\ge
1$ and $t\ge 0$. Then $\{h^N(\frac{k}{N}), k\ge 0 \}$ constitutes
the environment of the multiagent system. Thus we can define our
multiagent system which are driven by internal dynamics, interaction
among agents, and external economic fundamentals. Now we have
multiagent models evolving in random environment if we assume that
$\{h^N(\frac{k}{N}), k\ge 0 \}$ is a random sequence. The idea of
external economic fundamentals in deterministic or random
environment was used by Horst (2000).
\end{example}

We assume that for each $N\ge 1$, $A_N$ is an $\mathscr{L}_N$-valued
process defined on some probability space $(\Omega_N, \mathscr{F}_N,
Q_N)$, and $A$ is an $\mathscr{L}_c$-valued process defined on
$(\Omega, \mathscr{F}, Q)$. Then we can define multiagent model with
internal dynamics and random environment (MAMWIDARE), and multiagent
model with internal dynamics, multinomial sampling, and random
environment (MAMWIDMSARE).

Now, we state the joint and annealed weak convergence theorem for
the multiagent models evolving in random environment.
{\thm\label{thm:JointAnnealedMarginalRanCon} \textbf{(Joint and
annealed convergence)} Let $\{A_N\}$ and $A$ be defined above,
 $\mu_N\in \mathscr{P}(K_N)$ for each $N\ge 1$ and $\mu\in
 \mathscr{P}(K)$. Assume that $A_N\Rightarrow A$, and $\mu_N\Rightarrow
 \mu$.
\begin{itemize}
 \item[\textbf{1})] MAMWIDARE. For $N\ge 1$ and $\omega_N\in\Omega_N$,
 define $\Phi_N(A_N(\omega_N))\in \mathscr{P}(D_{K_N}[0,\infty))$ such that $\Phi_N(A_N(\omega_N))(Z_N(0))^{-1}=\mu_N$; for each $\omega\in \Omega$,
 define $\Phi(A(\omega))\in \mathscr{P}(D_{K}[0,\infty))$ such that $\Phi(A(\omega))(Z(0))^{-1}=\mu$. Then
\begin{align}
\label{eq:MAMWIDAREJointCong} (A_N, \Phi_N(A_N))\Rightarrow &(A, \Phi(A)) \\
\label{eq:MAMWIDAREAnnealCong} \lim_{N\to\infty}\rho(\int
\Phi(A_N(\omega_N))Q_N(d\omega_N), &
\int\Phi(A(\omega))Q(d\omega))=0.
\end{align}

\item[\textbf{2})] MAMWIDMSARE. For $N\ge 1$ and $\omega_N\in\Omega_N$,
 define $\Psi_N(A_N(\omega_N))\in \mathscr{P}(D_{K_N}[0,\infty))$ such that $\Psi_N(A_N(\omega_N))(Z_N(0))^{-1}=\mu_N$; for each $\omega\in \Omega$,
 define $\Psi(A(\omega))\in \mathscr{P}(D_{K}[0,\infty))$ such that $\Psi(A(\omega))(Z(0))^{-1}=\mu$. Then
\begin{align}
\label{eq:MAMWIDMSAREJointCong}
(A_N, \Psi_N(A_N))\Rightarrow &(A, \Psi(A)) \\
\label{eq:MAMWIDMSAREAnnealCong} \lim_{N\to\infty}\rho(\int
\Psi(A_N(\omega_N))Q_N(d\omega_N),
&\int\Psi(A(\omega))Q(d\omega))=0.
\end{align}
\end{itemize}

\section{Proof of Theorem \ref{thm:MAMIDTVE}}\label{se:ProofOfMAMWID}

At first, we consider some properties of
 $\{ {\mathbf n}^N(k), k \ge 0\}$ defined just by internal dynamics. Let $\mathbf{V}=(v_1,
\cdots, v_r)'$ be a positive
 vector. For any $k\ge 1$ and  $m\ge 0$, by (\ref{eq:XiNki}) and (\ref{eq:nk+In}), the independence of $\mathbf{\Xi}_{N,k,i}$'s,
  and the Markov property of $\{ {\mathbf n}^N(k), k \ge
0\}$, we have
\begin{equation}\label{eq:mgfforEvo}
E\biggl[\prod_{i=1}^{r}v_i^{n^N_i(m+k)}\biggl|{\mathbf
n}^N(m)\biggl]=\prod_{i=1}^{r}\biggl(P_{N,m,i\cdot}\prod_{l=1}^{k-1}P_{N,m+l}\mathbf{V}\biggl)^{n^N_i(m)},
\end{equation}
 where $\prod_{l=1}^{k-1}P_{N,m+l}=P_{N,m+1} \times\cdots \times
P_{N,m+k-1}$. We make the convention that when
 we denote the product of a sequence of matrices by prod, we actually make the multiplication from the left
  to the right as the index increases its order.

Then it follows by (\ref{eq:mgfforEvo}) that
\begin{equation}\label{eq:meanEvo}
E[{\mathbf n}^N(m+k)|{\mathbf n}^N(m)]={\mathbf
n}^N(m)\prod_{l=m}^{m+k-1}P_{N,l},
\end{equation} and for $1\le j\le
r$,
\begin{equation*}
\begin{aligned}
&\,E[n^N_j(m+k)(n^N_j(m+k)-1)|{\mathbf
n}^N(m)] \\
  =&\,\biggl({\mathbf n}^N(m)\biggl[\prod_{l=m}^{m+k-2}P_{N,l}\biggl]P_{N,m+k-1,\cdot,j}\biggl)^2 \\
  &\, -\sum_{i=1}^{r}\biggl( P_{N,m,i,\cdot}\biggl[\prod_{l=m+1}^{m+k-2}P_{N,l}\biggl]P_{N,m+k-1,\cdot,j}\biggl)^2 n^N_{i}(m),
\end{aligned}
\end{equation*}
where $P_{N,m+k-1,\cdot,j}$ is $j$-th column of matrix
$P_{N,m+k-1}$.
 Then we can get
\begin{equation}
\begin{aligned}\label{eq:squareOfDifEvo}
&\,E[(n^N_j(m+k)-n^N_j(m))^2|{\mathbf n}^N(m)]\\
  =&\,\biggl({\mathbf n}^N(m)\biggl[\prod_{l=m}^{m+k-2}P_{N,l}\biggl]P_{N,m+k-1,\cdot,j}\biggl)^2+{\mathbf n}^N(m)\biggl[\prod_{l=m}^{m+k-2}P_{N,l}\biggl]P_{N,m+k-1,\cdot,j} \\
  &\, -\sum_{i=1}^{r}\biggl( P_{N,m,i,\cdot}\biggl[\prod_{l=m+1}^{m+k-2}P_{N,l}\biggl]P_{N,m+k-1,\cdot,j}\biggl)^2 n^N_{i}(m) \\
  &\, -2{\mathbf n}^N(m)\biggl[\prod_{l=m}^{m+k-2}P_{N,l}\biggl]P_{N,m+k-1,\cdot,j}n^N_j(m)+(n^N_j(m))^2.
  \end{aligned}
  \end{equation}

For each $N\ge 1$, we define the transition operators on $\{\frac{
{\mathbf n}^N(k)}{N}:k\ge 0\}$ as follows:
\begin{equation}\label{eq:TranOpXN}
S_{N, k}f({\mathbf p})=E[f(\frac{ {\mathbf n}^N(k+1)}{N})|\frac{
{\mathbf n}^N(k)}{N}={\mathbf p}],
\end{equation} for each $k\ge 0$ and $f\in C(K_N)$,  ${\mathbf
p}\in K_N$.

{\lem\label{lem:SNNtGAtCon} Define $P_{N, k}$, and $\mathbf{X}^N$ in
subsection \ref{sse:MAMIDTVE}, and $S_{N,k}$ by (\ref{eq:TranOpXN}).
If $\lim_{N\to \infty}d_U(A_N, A)=0$, where $d_U$ is defined by
(\ref{eq:defndU}),
 then
\begin{equation}\label{eq:SNNtGAtCon}
 \lim_{N\to \infty} \sup_{0\le t\le T}\sup_{{\mathbf
p}\in K_N} |N[S_{N,[Nt]}-I]f({\mathbf p})-G_A(t)f({\mathbf p})|=0
\end{equation}
for any $f\in D= C^2(K)$. }

 \proof Let $f\in D= C^2(K)$, fix $T>0$. For $0\le t\le
T$, and ${\mathbf p}\in K_N$, note that by Taylor's expansion and
 (\ref{eq:meanEvo})
\begin{equation*}
\begin{aligned}
{}&[S_{N,[Nt]}-I]f({\mathbf p})\\
 =&\,E[(\mathbf{X}^N(t+\frac{1}{N})-{\mathbf p})\frac{\partial f}{\partial {\mathbf x}}({\mathbf p})'|\mathbf{X}^N(t)={\mathbf p}] \\
  &\, +E[(\mathbf{X}^N(t+\frac{1}{N})-{\mathbf p})\frac{\partial^2 f}{\partial {\mathbf x}^2}(\mathbf{X}^N_*(t))(\mathbf{X}^N(t+\frac{1}{N})-{\mathbf p})'|\mathbf{X}^N(t)={\mathbf p}] \\
 =&\,\frac{1}{N}{\mathbf p}A_N(\frac{[Nt]}{N})\frac{\partial f}{\partial {\mathbf x}}({\mathbf p})'\\
 &\, +E[(\mathbf{X}^N(t+\frac{1}{N})-{\mathbf p})\frac{\partial^2 f}{\partial {\mathbf x}^2}(\mathbf{X}^N_*(t))(\mathbf{X}^N(t+\frac{1}{N})-{\mathbf p})'|\mathbf{X}^N(t)={\mathbf p}],
\end{aligned}
\end{equation*}
where $\mathbf{X}^N_*(t)={\mathbf p}+\theta^N_t
(\mathbf{X}^N(t+\frac{1}{N})-{\mathbf p})$, for some $\theta^N_t\in
(0,1)$. Then
\begin{equation}\label{eq:NSNNtGAt}
\begin{aligned}
{}&|N[S_{N,[Nt]}-I]f({\mathbf p})-G_A(t)f({\mathbf p})| \\
 \le&\,|{\mathbf p}[A_N(\frac{[Nt]}{N})-A(t)]\frac{\partial f}{\partial {\mathbf x}}({\mathbf p})'| \\
 &\, +N|E[(\mathbf{X}^N(t+\frac{1}{N})-{\mathbf p})\frac{\partial^2 f}{\partial {\mathbf x}^2}(\mathbf{X}^N_*(t))(\mathbf{X}^N(t+\frac{1}{N})-{\mathbf p})'|\mathbf{X}^N(t)={\mathbf p}]|.
 \end{aligned}
 \end{equation}

Denote by $I_1(N,{\mathbf p},t)$ and $I_2(N,{\mathbf p},t)$, the
first and second term on the right hand side
 of (\ref{eq:NSNNtGAt}). Let $\|\frac{\partial^2 f}{\partial {\mathbf x}^2}\|=\max_{1\le i,j\le r}\|\frac{\partial^2 f}{\partial x_i \partial x_j }\|$. By H\"{o}lder inequality,
\begin{equation*}
\begin{aligned}
I_2(N,{\mathbf p},t)\le&\, \|\frac{\partial^2 f}{\partial {\mathbf x}^2}\| N \sum_{i=1}^r\sum_{j=1}^r E[|(X^N_i(t+\frac{1}{N})-p_i)(X^N_j(t+\frac{1}{N})-p_j)||\mathbf{X}^N(t)={\mathbf p}] \\
\le&\, \|\frac{\partial^2 f}{\partial {\mathbf x}^2}\| N
\sum_{i=1}^r\sum_{j=1}^r \biggl(E[(X^N_i(t+\frac{1}{N})-
p_i)^2|\mathbf{X}^N(t)={\mathbf p}]\\
 &\, \times E[(X^N_j(t+\frac{1}{N})-
p_j)^2|\mathbf{X}^N(t)={\mathbf p}]\biggl)^{\frac{1}{2}}.
\end{aligned}
\end{equation*}
 Since $A_N(t)$ is a constant on $[\frac{k}{N}, \frac{k+1}{N})$
for each $k\ge 0$, for fixed $1\le j\le r$, by
(\ref{eq:squareOfDifEvo}), we get
\begin{equation*}
\begin{aligned}
{}& N^2E[(X^N_j(t+\frac{1}{N})- p_j)^2|\mathbf{X}^N(t)={\mathbf p}] \\
  =&\, {\mathbf p}A_{N,\cdot,j}(\frac{[Nt]}{N})+({\mathbf p}A_{N,\cdot,j}(\frac{[Nt]}{N}))^2-2p_j a_{N,j,j}(\frac{[Nt]}{N})- \sum_{k=1}^r \frac{p_k}{N} a_{N,k,j}^2 (\frac{[Nt]}{N}) \\
 \le &\, \sum_{l=1}^r p_l |A_{N,l,j}(t)|+(\sum_{l=1}^r p_l |A_{N,l,j}(t)|)^2+2p_j |a_{N,j,j}(t)|\\
 \le &\, \sum_{l=1}^r  |A_{N,l,j}(t)|+(\sum_{l=1}^r  |A_{N,l,j}(t)|)^2+2 |a_{N,j,j}(t)|.
 \end{aligned}
 \end{equation*}

Since $\lim_{N\to \infty}d_U(A_N, A)=0$ and $A$ is bounded on
$[0,T]$, there exists $C_T>0$ such that
\begin{equation*}
N^2E[(X^N_j(t+\frac{1}{N})- p_j)^2|\mathbf{X}^N(t)={\mathbf p}]\le
C_T
\end{equation*}
 for any $N\ge 1$, $0\le t\le T$, ${\mathbf p}\in
K_N$, and $1\le j\le r$. Then
\begin{equation}\label{I2NPtBd}
\begin{aligned}
\sup_{0\le t\le T}\sup_{{\mathbf p}\in K_N}I_2(N,{\mathbf
p},t)\le&\, \frac{1}{N}\|\frac{\partial^2 f}{\partial {\mathbf
x}^2}\|r^2 C_T.
\end{aligned}
\end{equation}
 Thus $\lim_{N\to
\infty} \sup_{0\le t\le T}\sup_{{\mathbf p}\in K_N}I_2(N,{\mathbf
p},t)=0$.

For matrix $B=(b_{i,j})_{r\times r}$ and vector ${\mathbf
y}=(y_1,\cdots,y_r)'$, let $\|B\|=\sum_{i,j=1}^r|b_{i,j}|$ and
$\|{\mathbf y}\|$ be the Euclidean norm of $\mathbf{y}$. Then $\|B
{\mathbf y}\|\le \|B\| \cdot\|\mathbf{y}\|$. It follows that
\begin{equation}\label{eq:I1Npt}
\begin{aligned}
I_1(N,{\mathbf p},t)\le&\, \|{\mathbf
p}'\|\cdot\|A_N(\frac{[Nt]}{N})-A(t)\|\cdot
\|\frac{\partial f}{\partial {\mathbf x}}({\mathbf p})'\|\\
\le &\,\|A_N(\frac{[Nt]}{N})-A(t)\|\cdot \|\frac{\partial
f}{\partial {\mathbf x}}({\mathbf p})'\|.
\end{aligned}
\end{equation}
 Notice that
$A_N(\frac{[Nt]}{N})=A_N(t)$, $t\ge 0$, thus
\begin{equation}\label{IN1tGenearlMAM}
\lim_{N\to \infty} \sup_{0\le t\le
T}\sup_{{\mathbf p}\in K_N}I_1(N,{\mathbf p},t)=0
\end{equation}
 follows by (\ref{eq:I1Npt}) and $\lim_{N\to \infty}d_U(A_N, A)=0$. Therefore
(\ref{eq:SNNtGAtCon}) is proved.

{\cor\label{lem:NSNNtBd&SNNtGoTo0} With the same conditions as those
in Lemma \ref{lem:SNNtGAtCon}, we have
\begin{equation}\label{eq:NSNNtBd}
\sup_N \sup_{0\le t \le
T}\sup_{{\mathbf p}\in K_N} |N[S_{N,[Nt]}-I]f({\mathbf p})|<\infty
\end{equation}
and
\begin{equation}\label{eq:SNNtGoTo0}
\lim_{N\to \infty} \sup_{0\le t \le T}\sup_{{\mathbf p}\in K_N}
|[S_{N,[Nt]}-I]f({\mathbf p})|=0
\end{equation}
 for any $T>0$ and $f\in C^2(K)$. }

\proof This follows from Lemma \ref{lem:SNNtGAtCon} immediately.

{\rem\label{rem:weakerCond} We can use Taylor's expansion of order 1
to prove directly that (\ref{eq:NSNNtBd}) and (\ref{eq:SNNtGoTo0})
hold for any $T>0$ and $f\in C^1(K)$.}

{\lem\label{lem:UniGInhomWOMS} Define the time-dependent generator
$\{G_A(t), 0\le t<\infty\} $ on $C^1(K)$ by
(\ref{eq:genGeneralInhom}). Let $\mu\in \mathscr{P}(K)$. Then the
$D_K[0,\infty)$ martingale problem for $(G_A, \mu)$ has at most one
solution. }

\proof The limit process $\mathbf{X}^{\infty}$ of Theorem
\ref{thm:MAMIDTVE} is deterministic such that
$(\mathbf{X}^{\infty})'$ satisfies the linear differential equation
$$ \frac{d \mathbf{x}(t)}{d t}=A'(t)\mathbf{x}(t).$$
Then we have the uniqueness of the $D_K[0,\infty)$ martingale
problem for $(G_A, \mu)$.

{ {\bf Proof of Theorem \ref{thm:MAMIDTVE}, Part \textbf{1}).}} This
part follows by Lemma \ref{lem:SNNtGAtCon}, \ref{lem:UniGInhomWOMS},
and Remark \ref{rem:GenAndSimp}.

Next we make preparations for proving Theorem \ref{thm:MAMIDTVE},
Part \textbf{2}).

At first, we consider the multiagent model with only multinomial
sampling (MAMWMS). This is illustrated by Figure \ref{MultSamp:fig}.
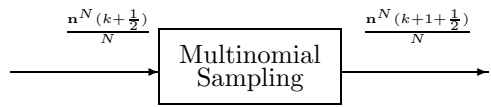
\begin{figure}[h]
\begin{center}
\vspace{1.5cm}
\begin{picture}(400,0)
$\stackrel{\hspace{4.0cm}\frac{\mathbf{n}^N(k+\frac{1}{2})}{N}}{\put(20.0,12){\vector(1,0){56.5}}}$
\put(3,0){\framebox(68.5,28)[c]{$\stackrel{\mbox{Multinomial}}{\mbox{Sampling}}$}
}$\stackrel{\hspace{2.8cm}\frac{\mathbf{n}^N(k+1+\frac{1}{2})}{N}}{\put(11.0,12){\vector(1,0){55.5}}}$
\end{picture}

\end{center}
\caption{The transition structure of MAMWMS} \label{MultSamp:fig}
\end{figure}

This model is similar to the neutral Wright-Fisher model of
population genetics. However in contrast to the standard genetics
model here the internal state change (which would correspond to
mutation) can be influenced by an external random environment. If we
don't include in the Wright-Fisher model {\em selection} and {\em
gene mutation}, we get a model which is very close to MAMWMS. We
define a transition operator $T_N$ related to the homogeneous Markov
chain $\{\frac{\mathbf{ n}^N(k+\frac{1}{2})}{N}, k=0,1, \cdots\}$ on
$C(K_N)$ as follows:
\begin{equation*}
T_N f({\mathbf p})=E[f(\frac{\mathbf{
n}^N(1+\frac{1}{2})}{N})|\frac{{\mathbf
n}^N(\frac{1}{2})}{N}={\mathbf p}], \mbox{ }f\in C(K_N), \mbox{
}{\mathbf p}\in K_N.
\end{equation*}
Let
\begin{equation*}
\mathbf{Z}^N(t)=\frac{\mathbf{ n}^N([Nt]+\frac{1}{2})}{N}
\end{equation*} and let $\mathbf{Z}$ be a diffusion process in $K$ with the generator
$G_B$ defined by (\ref{eq:DefnOfB}). Then we have the following
analog of the classical Wright-Fisher diffusion limit.
{\prop With
the conditions above,
\begin{equation}\label{eq:ConvTN}
\lim_{N\to \infty} \sup_{{\mathbf p}\in K_N} |N(T_N-I)f({\mathbf
p})-G_Bf({\mathbf p})|=0
\end{equation} for every $f\in C^2(K)$. If
$\mathbf{Z}^N(0)\Rightarrow \mathbf{Z}(0)$ in $K$, then
$\mathbf{Z}^N\Rightarrow \mathbf{Z}$ in $D_K[0,\infty)$.}

Define $\{U_{N,k}, k\ge 0\}$ as follows:
\begin{equation*}
U_{N,k} f({\mathbf p})=E[f(\frac{{\mathbf
n}^N(k+1)}{N})|\frac{{\mathbf n}^N(k)}{N}={\mathbf p}], \mbox{ }f\in
C(K_N), \mbox{ }{\mathbf p}\in K_N \mbox{ for any $k\ge 0$}.
\end{equation*}
 Then by the Figure \ref{MAMWIDAMS:fig}, since $S_{N,k}$ and $T_N$ are one step transition operators related to the internal
 dynamics and the global multinomial sampling respectively, we have
\begin{equation}\label{eq:UTSNk}
 U_{N,k}= S_{N,k} T_N.
\end{equation}

 {\lem\label{lem:ConvUNK} Assume that $\lim_{N\to
\infty}d_U(A_N, A)=0$. With the definitions above, we have
\begin{equation}\label{eq:ConvUNk}
\lim_{N\to \infty} \sup_{{\mathbf p}\in K_N}
|N[U_{N,[Nt]}-I]f({\mathbf p})-G_{AB}(t)f({\mathbf p})|=0
\end{equation} for every $f\in C^3(K)$ and $t\ge 0$.
}

\proof Let $f\in C^3(K)$, and fix $t\ge 0$. Choose $T$ such that
$T\ge t$, and ${\mathbf p}\in K_N$, notice that by (\ref{eq:UTSNk})
\begin{equation*}
\begin{aligned}{}& |N[U_{N, [Nt]}-I]f({\mathbf p})-G_{AB}(t)f({\mathbf p})|  \\
   =&\,|[NS_{N, [Nt]}(T_N-I)]f({\mathbf p}) - G_B f({\mathbf p}) + N[S_{N, [Nt]}-I]f({\mathbf p}) -G_A(t) f({\mathbf p})| \\
\le &\, |S_{N, [Nt]}[N(T_N-I)]f({\mathbf p}) - S_{N, [Nt]} G_B
f({\mathbf p})|+| S_{N, [Nt]} G_B f({\mathbf p})- G_B f({\mathbf
p})|\\
   &\,+|N[S_{N, [Nt]}-I]f({\mathbf p}) -G_A(t) f({\mathbf p})|.
 \end{aligned}
 \end{equation*}
 It follows by $S_{N, [Nt]}$ being a contraction that
\begin{equation*}
\begin{aligned}{}& \sup_{{\mathbf p}\in K_N}
|N[U_{N,[Nt]}-I]f({\mathbf p})-G_{AB}(t)f({\mathbf p})| \\
\le &\, \sup_{{\mathbf p}\in K_N} |[N(T_N-I)]f({\mathbf p})-G_B
f({\mathbf p})|+\sup_{{\mathbf p}\in K_N} | S_{N, [Nt]} G_B
f({\mathbf p})- G_B f({\mathbf
p})|\\
   &\,+\sup_{{\mathbf p}\in K_N} |N[S_{N, [Nt]}-I]f({\mathbf p}) -G_A(t) f({\mathbf p})|.
\end{aligned}
 \end{equation*}

(\ref{eq:ConvUNk}) is then proved by (\ref{eq:ConvTN}),
(\ref{eq:SNNtGoTo0}), Remark \ref{rem:weakerCond}, and
(\ref{eq:SNNtGAtCon}).

{\lem\label{lem:NUNNtBd} With the same conditions as those in Lemma
\ref{lem:ConvUNK}, we have
\begin{equation}\label{eq:NUNNtBd}
\sup_N\sup_{0\le t \le T}\sup_{{\mathbf p}\in K_N}
|N[U_{N,[Nt]}-I]f({\mathbf p})|<\infty
\end{equation}
 for any $T>0$ and $f\in C^3(K)$. }

\proof Let $f\in C^3(K)$, fix $T>0$. For $0\le t\le T$, and
${\mathbf p}\in K_N$, notice that
\begin{equation*}
\begin{aligned}\sup_{{\mathbf p}\in K_N}|N[U_{N, [Nt]}-I]f({\mathbf p})|=&\,\sup_{{\mathbf p}\in K_N}|[NS_{N, [Nt]}(T_N-I)]f({\mathbf p})+N[S_{N, [Nt]}-I]f({\mathbf p})|  \\
\le &\, \sup_{{\mathbf p}\in K_N}|N[T_N-I]f({\mathbf
p})|+\sup_{{\mathbf p}\in K_N}|N[S_{N, [Nt]}-I]f({\mathbf p})|.
\end{aligned}
\end{equation*}
(\ref{eq:NUNNtBd}) follows by (\ref{eq:NSNNtBd}) and
(\ref{eq:ConvTN}).

Now we introduce some notations before we prove the uniqueness of
the martingale problem for $(G_{AB}, \mu)$, $\mu\in \mathscr{P}(K)$
on $C^3(K)$. For each $n\ge 1$, let $\mathbb{N}_n=\{\mathbf{\alpha}:
\mathbf{\alpha}\in (Z_+)^r, \sum_{i=1}^r \alpha_i=n\}$. It is clear
that $\mathbb{N}_n$ has $c_n=\bin{n+r-1}{\hspace{4.5mm}n}$ elements.
Define an order $'>'$ on $\mathbb{N}_n$ as follows: let $\alpha$,
$\alpha' \in \mathbb{N}_n$, $\alpha>\alpha'$ $\Leftrightarrow$
$\alpha_1>\alpha'_1$ or there exists $2\le k\le r$, such that
$\alpha_l=\alpha'_l$ for $1\le l\le k-1$ and $\alpha_k>\alpha'_k$.
We use $\mathbb{N}_n$ as the index set when we define
$D_{R^{c_n\times c_n}}[0,\infty)$-valued matrix functions. For given
$A\in D_{R^{r\times r}}[0,\infty)$, define
$A^{(n)}=(a^{(n)}_{\alpha,\alpha'})_{\alpha,\alpha'\in \mathbb{N}_n}
\in D_{R^{c_n\times c_n}}[0,\infty)$ as follows: for each $\alpha$,
$\alpha'\in \mathbb{N}_n$ and $t\ge 0$,
\begin{equation}\label{eq:defnofAn}
a^{(n)}_{\alpha,\alpha'}(t)=\left\{
\begin{array}{ll}
        \sum_{j=1}^r \alpha_j a_{jj}(t), &\mbox{if } \alpha=\alpha'; \\
       \alpha_j a_{kj}(t),     &\mbox{if there exist $1\le j, k\le r$, $j\not=k$, such that } \alpha_l=\alpha'_l,  \\
       &\mbox{ for } 1\le l\le r \mbox{ satisfying $l\not=j$ or $k$}, \mbox{ and } \alpha_j'=\alpha_j-1,   \\
       &\alpha'_k=\alpha_k+1 ; \\
                     0, &\mbox{otherwise}.
 \end{array} \right.
\end{equation}
 We arrange the elements of $A^{(n)}$ along the rows and columns
decreasingly by the order $'>'$. For each $n\ge 1$, we define a
$c_n\times c_n$ matrix
$B^{(n)}=(b^{(n)}_{\alpha,\alpha'})_{\alpha,\alpha'\in
\mathbb{N}_n}$ as follows: for each $\alpha$, $\alpha'\in
\mathbb{N}_n$,
\begin{equation}\label{eq:defnofBn}
b^{(n)}_{\alpha,\alpha'}=\left\{
\begin{array}{ll}
                               \frac{1}{2}(n-n^2), &\mbox{if } \alpha=\alpha'; \\
                                                0, &\mbox{otherwise}.
 \end{array} \right.
\end{equation}
For each $n\ge 2$,  using $\mathbb{N}_n$ and $\mathbb{N}_{n-1}$ as
the index sets, we define a $c_n\times c_{n-1}$ matrix
$D^{(n)}=(d^{(n)}_{\alpha,\alpha'})_{\alpha\in
\mathbb{N}_n,\alpha'\in \mathbb{N}_{n-1}}$ as follows: for each
$\alpha\in \mathbb{N}_n$, $\alpha'\in \mathbb{N}_{n-1}$,
\begin{equation}\label{eq:defnofdn}
d^{(n)}_{\alpha,\alpha'}=\left\{ \begin{array}{ll}
                                                     \frac{1}{2}\alpha_i(\alpha_i-1), &\mbox{if $\alpha'_i=\alpha_i-1$ and $\alpha'_l=\alpha_l$ for $l\not=i$, $1\le l\le r$}; \\
                                                       0, &\mbox{otherwise}.
 \end{array} \right.
\end{equation}
We also arrange the elements of $B^{(n)}$, $D^{(n)}$ along the rows
and columns decreasingly by the order $'>'$.

{\lem\label{lem:UniGInhomWithMS} Define the time-dependent generator
$\{G_{AB}(t), 0\le t<\infty\} $ on $C^2(K)$ by
(\ref{eq:defnofGABs}). Let $\mu\in \mathscr{P}(K)$. Then the
$D_K[0,\infty)$ martingale problem for $(G_{AB}, \mu)$ on $C^3(K)$
has at most one solution. }

\proof Assume that $\{\mathbf{Y}(t)=(Y_1(t), \cdots, Y_r(t)):0\le
t<\infty\}$ is one solution of the martingale problem for $(G_{AB},
\mu)$ on $C^3(K)$. Let $n\ge 1$, for arbitrary $\alpha\in
\mathbb{N}_n$, and $t\ge 0$, define $f_{\alpha}(x_1, \cdots,
x_r)=x_1^{\alpha_1}\cdots x_r^{\alpha_r}$ and
$y^{\alpha}_n(t)=E[f_{\alpha}(\mathbf{Y}(t))]$. Define the column
vector $\mathbf{y}_n(t)=(y^{\alpha}_n(t))_{\alpha\in \mathbb{N}_n}$,
where the elements of $\mathbf{y}_n(t)$ is arranged decreasingly by
the order $'>'$.

Then, for $n\ge 1$ and given $\alpha\in \mathbb{N}_n$ and $t>0$, we
have
\begin{equation}\label{eqGABFirst}
E[f_{\alpha}(\mathbf{Y}(t))]=E[f_{\alpha}(\mathbf{Y}(0))]+\int_0^t
E[G_{AB}(s)f_{\alpha}(\mathbf{Y}(s))]ds.
\end{equation}

At first, for $n=1$, since $G_B f_{\alpha}\equiv 0$ for any
$\alpha\in \mathbb{N}_1$, we have
\begin{equation*}
E[f_{\alpha}(\mathbf{Y}(t))]=E[f_{\alpha}(\mathbf{Y}(0))]+\int_0^t
E[G_A(s)f_{\alpha}(\mathbf{Y}(s))]ds,
\end{equation*}
which implies
\begin{equation}\label{eq:alphaeqnequa1}
\mathbf{y}_1(t)=\mathbf{y}_1(0)+\int_0^t
A^{(1)}(s)\mathbf{y}_1(s)ds.
\end{equation}

Next, we calculate $E[f_{\alpha}(\mathbf{Y}(t))]$ for $n\ge 2$. By
(\ref{eq:DefnOfB}), (\ref{eq:defnofBn}) and (\ref{eq:defnofdn}), we
can get
\begin{equation}\label{eq:GBfWMS}
\begin{aligned}
&\, E[G_Bf_{\alpha}(\mathbf{Y}(s))]\\
=&\,\frac{1}{2}E[ \sum_{i=1}^r b_{ii}(\mathbf{Y}(s))\frac{\partial^2 f_{\alpha}}{\partial x_i^2}(\mathbf{Y}(s))]+\frac{1}{2}E[ \sum_{i=1}^r\sum_{j=1,j\not=i}^r b_{ij}(\mathbf{Y}(s))\frac{\partial^2 f_{\alpha}}{\partial x_i x_j}(\mathbf{Y}(s))] \\
=&\,\frac{1}{2}\sum_{i=1}^r E[ Y_i(s)(1-Y_i(s))\alpha_i(\alpha_i-1)Y_1(s)^{\alpha_1}\cdots Y_i(s)^{\alpha_i-2}\cdots Y_r(s)^{\alpha_r}] \\
&\, -\frac{1}{2} \sum_{i=1}^r\sum_{j=1, j\not=i}^r E[ \alpha_i\alpha_j Y_1(s)^{\alpha_1}\cdots Y_r(s)^{\alpha_r}] \\
=&\,\frac{1}{2}\sum_{i=1}^r \alpha_i(\alpha_i-1)E[Y_1(s)^{\alpha_1}\cdots Y_i(s)^{\alpha_i-1}\cdots Y_r(s)^{\alpha_r}]\\
&\, -\frac{1}{2}[\sum_{i=1}^r \alpha_i(\alpha_i-1)+  \sum_{i=1}^r\sum_{j=1, j\not=i}^r  \alpha_i\alpha_j ]E[ Y_1(s)^{\alpha_1}\cdots Y_r(s)^{\alpha_r}] \\
=&\, D^{(n)}_{\alpha,\cdot}\mathbf{y}_{n-1}(s)-\frac{1}{2}[\sum_{i=1}^r\sum_{j=1}^r\alpha_i\alpha_j-\sum_{i=1}^r \alpha_i] E[ Y_1(s)^{\alpha_1}\cdots Y_r(s)^{\alpha_r}] \\
=&\, D^{(n)}_{\alpha,\cdot}\mathbf{y}_{n-1}(s)-\frac{1}{2}(n^2-n) E[ Y_1(s)^{\alpha_1}\cdots Y_r(s)^{\alpha_r}]  \\
=&\, D^{(n)}_{\alpha,\cdot}\mathbf{y}_{n-1}(s)+
B^{(n)}_{\alpha,\cdot}\mathbf{y}_{n}(s),
\end{aligned}
\end{equation}
 where $D^{(n)}_{\alpha,\cdot}$ is the $\alpha$-row of the matrix
$D^{(n)}$, $B^{(n)}_{\alpha,\cdot}$ the $\alpha$-row of the matrix
$B^{(n)}$. By (\ref{eq:defnofAn}) we can get
\begin{equation}\label{eq:alphatheq}
\begin{aligned}
\int_0^t
E[G_A(s)f_{\alpha}(\mathbf{Y}(s))]ds=&\,\int_0^t E[\mathbf{Y}(s) A(s)\frac{\partial f_{\alpha}}{\partial \mathbf{x}} (\mathbf{Y}(s))']ds \\
=&\,\int_0^t \biggl[ \sum_{j=1}^r \alpha_j a_{j,j}(s) y^{\alpha}_n(s)+\sum_{k=1}^r\sum_{j=1}^{k-1} \alpha_j a_{k,j}(s) y^{(\alpha_1,\cdots,\alpha_j-1,\cdots,\alpha_k+1,\cdots,\alpha_r)}(s) \\
&\, +\sum_{k=1}^r\sum_{j=k+1}^{r} \alpha_j a_{k,j}(s) y^{(\alpha_1,\cdots,\alpha_k+1,\cdots,\alpha_j-1,\cdots,\alpha_r)}(s)\biggl]ds \\
=&\,\int_0^t A^{(n)}_{\alpha,\cdot}(s)\mathbf{y}_n(s)ds,
\end{aligned}
\end{equation}
 where $A^{(n)}_{\alpha,\cdot}$ is the $\alpha$-row of the matrix
$A^{(n)}$. Then by (\ref{eqGABFirst}), (\ref{eq:GBfWMS}), and
(\ref{eq:alphatheq}), for $n\ge 2$, we obtain
\begin{equation}\label{eq:alphaeqAnBnDn}
\mathbf{y}_n(t)=\mathbf{y}_n(0)+\int_0^t
D^{(n)}\mathbf{y}_{n-1}(s)ds+\int_0^t
[A^{(n)}(s)+B^{(n)}]\mathbf{y}_n(s)ds.
\end{equation}

Define for $n=1$,
\begin{equation*}
V_1(s,t)= \left\{
\begin{array}{ll}
                                    A^{(1)}(t)  & \mbox{if $0\le t\le s<\infty$}, \\
                                     0,         & \mbox{otherwise};
 \end{array} \right.
 \end{equation*}
and for $n\ge 2$,
\begin{equation*}
V_n(s,t)= \left\{
\begin{array}{ll}
                                    A^{(n)}(t)+B^{(n)}  & \mbox{if $0\le t\le s<\infty$}, \\
                                     0,         & \mbox{otherwise}.
 \end{array} \right.
 \end{equation*}
Let $f_1(t)=\mathbf{y}_1(0)$, $0\le t<\infty$. Then, by
(\ref{eq:alphaeqnequa1})  and the theory of {\em Volterra equations
of the second kind} (see e.g. Smithies (1958) or Tricomi (1957)), we
know that $\mathbf{x}_1=\mathbf{y}_1$ is the unique solution of the
{\em Volterra equation} of the second kind
\begin{equation}\label{eq:VoltEqVneq1}
\mathbf{x}(s)=f_1(s)+\int_0^s V_1(s,t)\mathbf{x}(t)dt, \,\,(0\le
s\le T)
\end{equation}
 on the space $\mathbb{L}^2([0,T], R^{c_1})$
for any $T>0$. We can do this procedure recursively. Assume that we
know that for $n\ge 1$ the unique solution
$\mathbf{x}_n=\mathbf{y}_n$ is determined, then we define
$f_{n+1}(t)=\mathbf{y}_{n+1}(0)+\int_0^t
D^{(n+1)}\mathbf{x}_{n}(s)ds$. Then by (\ref{eq:alphaeqAnBnDn}), we
know that $\mathbf{x}_{n+1}=\mathbf{y}_{n+1}$ is the unique solution
of the {\em Volterra equation} of the second kind
\begin{equation*}
 \mathbf{x}(s)=f_{n+1}(s)+\int_0^s
V_{n+1}(s,t)\mathbf{x}(t)dt, \,\,(0\le s\le T)
\end{equation*}
 on the space $\mathbb{L}^2([0,T], R^{c_{n+1}})$ for any $T>0$. Note that
the construction of the {\em Volterra equations} does depend just on
$\mathbf{y}_n(0)$ for $n\ge 1$, and does not depend on
$\mathbf{y}_n(t)$ for $t>0$, $n\ge 1$. Then we conclude that all
moments of the one-dimensional marginal distribution of any two
solutions of the martingale problem for $(G_{AB}, \mu)$ are the
same. The uniqueness of the martingale problem for $(G_{AB}, \mu)$
on $C^3(K)$ then follows by Theorem 4.2, Chapter 4 of Either and
Kurtz (1986).

{ {\bf Proof of Theorem \ref{thm:MAMIDTVE}, Part \textbf{2)}.} This
part follows by Lemma \ref{lem:ConvUNK}, \ref{lem:NUNNtBd},
\ref{lem:UniGInhomWithMS} and Corollary \ref{InhomMarkovChainCon}.

\section{Multiagent models in random environment}\label{se:MAMIRE}

\subsection{Measurability with respect to random
environment}\label{sse:mwrtre}

At first, for $\Phi$ and $\Psi$ define in Subsection
\ref{sse:MAMSIRE}, we have the following Lemma.

{\lem\label{thm:ContiPhiPsi}
  $\Phi$ and $\Psi$ are continuous mappings from $(\mathscr{L}_c, d_U)$ to
$(\mathscr{P}(D_K[0,\infty)),\rho)$. }

\proof This is immediate by
 (\ref{eq:genGeneralInhom}), (\ref{eq:defnofGABs}) and Corollary \ref{InhomMarkovGntGtCon}.

Now we consider the measurability related to $A$ if $A$ is an
$\mathscr{L}_c$-valued process.

{\lem\label{lem:measA} Assume that
$A$ is an $\mathscr{L}_c$-valued process defined on some probability
space
 $(\Omega, \mathscr{F}, Q)$. Then $(A,\Phi(A))$, $(A,\Psi(A))$ are $\mathscr{F}/[\mathscr{B}(\mathscr{L}_c)\otimes\mathscr{B}(\mathscr{P}(D_K[0,\infty)))]$-measurable. }

\proof Since $K$ is separable, so are $D_K[0,\infty)$ and
$\mathscr{P}(D_K[0,\infty))$. Since $\mathscr{L}_c$ is also
 separable, to prove that $(A,\Phi(A))$ is measurable, it suffices to prove that if $C\in \mathscr{B}(\mathscr{L}_c)$,
  $D\in \mathscr{B}(\mathscr{P}(D_K[0,\infty)))$, $(A,\Phi(A))^{-1}(C\times D)\in \mathscr{F}$. This is clear by
  Lemma \ref{thm:ContiPhiPsi}, and $(A,\Phi(A))^{-1}(C\times D)=A^{-1}(C\cap \Phi^{-1}(D))$. Similarly, we can
    prove that $(A,\Psi(A))$ is measurable.

Next, we consider the measurability related to $\{A_N\}$. Since
$\{A_N(\frac{k}{N}), k\ge 0\}$
          determines the transition structure of a Markov chain by $\{S_{N,k},k\ge 0\}$ in MAMWID or by
           $\{U_{N,k},k\ge 0\}$ in MAMWIDAMS, it is clear that $\Phi_N$ and $\Psi_N$ are continuous mappings
            from $\mathscr{L}_N$ to $\mathscr{P}(D_{K_N}[0,\infty))$. Then we
            have the following lemma.

{\lem\label{lem:measAN} Assume that $A_N$ is an
$\mathscr{L}_N$-valued process defined on some probability space
 $(\Omega_N, \mathscr{F}_N, Q_N)$. Then $(A_N,\Phi_N(A_N))$, and $(A_N,\Psi_N(A_N))$ are
  $\mathscr{F}/[\mathscr{B}(\mathscr{L}_N)\otimes\mathscr{B}(\mathscr{P}(D_{K_N}[0,\infty)))]$-measurable. }

{\rem\label{rem:ReintprOfConvGen} Using the notations $\Phi_N$,
$\Psi_N$, $\Phi$, and $\Psi$, we can restate Theorem
\ref{thm:MAMIDTVE} as follows. Let $\hat{A}_N\in \mathscr{L}_N$, and
  $\hat{A}\in \mathscr{L}_c$. Let $\mu_N\in \mathscr{P}(K_N)$, and $\mu\in
  \mathscr{P}(K)$. Assume that $\lim_{N\to\infty}d(\hat{A}_N, \hat{A})=0$, and $\mu_N \Rightarrow
  \mu$.
\begin{itemize}
\item[\textbf{1})] MAMWID. Define $\Phi_N(\hat{A}_N)\in
\mathscr{P}(D_{K_N}[0,\infty))$ and $\Phi(\hat{A})\in
\mathscr{P}(D_{K}[0,\infty))$, such that
$\Phi_N(\hat{A}_N)Z_N(0)^{-1}=\mu_N$ and
 $\Phi(\hat{A})Z(0)^{-1}=\mu$.  Then \newline
  $\lim_{N\to\infty}\rho(\Phi_N(\hat{A}_N), \Phi(\hat{A}))=0$.

\item[\textbf{2})] MAMWIDAMS. Define $\Psi_N(\hat{A}_N)\in
\mathscr{P}(D_{K_N}[0,\infty))$ and $\Psi(\hat{A})\in
\mathscr{P}(D_{K}[0,\infty))$, such that
$\Psi_N(\hat{A}_N)Z_N(0)^{-1}=\mu_N$ and
$\Psi(\hat{A})Z(0)^{-1}=\mu$. Then \newline
 $\lim_{N\to\infty}\rho(\Psi_N(\hat{A}_N), \Psi(\hat{A}))=0$.
\end{itemize} }

\subsection{Proof of Theorem \ref{thm:JointAnnealedMarginalRanCon}}\label{sse:POTJA}

The next two lemmas will be used in the proof of Theorem
\ref{thm:JointAnnealedMarginalRanCon}.

{\lem\label{lem:condi12} Let $\hat{A}=(\hat{a}_{i,j})_{r\times r}$
be an $\mathscr{L}_c$-valued process a.s., defined on some
probability space $(\hat{\Omega}, \hat{\mathcal{F}}, \hat{P})$, and
$\tilde{A}=(\tilde{a}_{i,j})_{r\times r}$ be a $C_{R^{r\times
r}}[0,\infty)$-valued process defined on some probability space
$(\tilde{\Omega}, \tilde{\mathscr{F}}, \tilde{P})$. If
$\hat{P}\hat{A}^{-1}=\tilde{P}\tilde{A}^{-1}$, then $\tilde{A}$ is
also an $\mathscr{L}_c$-valued process a.s. }

\proof The proof of this lemma is straightforward.

{\lem\label{lem:stepwiseConst} Let $A_N$ be an
$\mathscr{L}_N$-valued process defined on some probability space
\newline $(\Omega_N, \mathscr{F}_N, Q_N)$ a.s., and let $\hat{A}_N$ be a
$D_{R^{r\times r}}[0,\infty)$-valued process defined on
$(\hat{\Omega}_N, \hat{\mathscr{F}}_N, \hat{Q}_N)$. If $Q_N
A_N^{-1}=\hat{Q}_N \hat{A}_N^{-1}$, then $\hat{A}_N$ is also an
$\mathscr{L}_N$-valued process a.s. }

\proof We omit this proof.

Lemma \ref{lem:condi12} and \ref{lem:stepwiseConst} indicate that
the conditions specified at the beginning of subsection
\ref{sse:MAMIDTVE} for the processes $\{A_N\}$ and $A$ just depend
on the distributions of $\{A_N\}$ and $A$.

{ {\bf Proof of Theorem \ref{thm:JointAnnealedMarginalRanCon}.} We
just prove part \textbf{1)}. Part \textbf{2)} can be proved
similarly.

At first, we prove (\ref{eq:MAMWIDAREJointCong}). Since $\{A_N\}$
are $D_{R^{r\times r}}[0,\infty)$-valued processes, $A$ is a
$C_{R^{r\times r}}[0,\infty)$-valued process, and $A_N\Rightarrow
A$, by Skorohod Representation theorem, there exist some probability
space $(\hat{\Omega}, \hat{\mathscr{F}}, \hat{Q})$ and a sequence of
 $D_{R^{r\times r}}[0,\infty)$-valued processes $\{\hat{A}_N\}$ and a $C_{R^{r\times r}}[0,\infty)$-valued process $\hat{A}$ satisfying
  $\lim_{N\to\infty} d(\hat{A}_N,\hat{A})=0$ $\hat{Q}$-a.s, $\hat{Q}\hat{A}^{-1}=Q A^{-1}$, and
   $\hat{Q}\hat{A}_N^{-1}=Q_N A_N^{-1}$ on $\mathscr{B}(D_{R^{r\times r}}[0,\infty))$ for each $N\ge 1$.
    By Lemma \ref{lem:condi12} and \ref{lem:stepwiseConst}, $\hat{A}$ is an $\mathscr{L}_c$-valued process,
    and $\hat{A}_N$ is an $\mathscr{L}_N$-valued process for each $N\ge 1$. By Lemma \ref{lem:measA} and
    \ref{lem:measAN}, $(\hat{A}, \Phi(\hat{A}))$ and $(\hat{A}_N, \Phi_N(\hat{A}_N))$ are measurable.
     By Remark \ref{rem:ReintprOfConvGen}, we have
\begin{equation}\label{eq:PhiANPHiA}
\lim_{N\to\infty}\rho(\Phi_N(\hat{A}_N(\hat{\omega})),
\Phi(\hat{A}(\hat{\omega})))=0, \mbox{ $\hat{Q}$-a.s.}
\end{equation}

For each $f\in \bar{C}(\mathscr{L}_c)$, and $g\in
\bar{C}(\mathscr{P}(D_K[0,\infty)))$, by bounded convergence theorem
\begin{equation*}
\begin{aligned}
&\,\lim_{N\to\infty}\int f(A_N(\omega_N))g(\Phi_N(A_N(\omega_N)))Q(d\omega_N)\\
=&\,\lim_{N\to\infty}\int f(\hat{A}_N(\hat{\omega}))g(\Phi_N(\hat{A}_N(\hat{\omega})))\hat{Q}(d\hat{\omega}) \\
 =&\,\int f(\hat{A}(\hat{\omega}))g(\Phi(\hat{A}(\hat{\omega})))\hat{Q}(d\hat{\omega})  \\
 =&\,\int f(A(\omega))g(\Phi(A(\omega)))Q(d\omega).
 \end{aligned}
\end{equation*}
Since $\mathscr{L}_c$ and $\mathscr{P}(D_K[0,\infty))$ are
separable, $\bar{C}(\mathscr{L}_c)$ and
$\bar{C}(\mathscr{P}(D_K[0,\infty)))$
 are convergence determining on $(\mathscr{L}_c, d_U)$ and $(\mathscr{P}(D_K[0,\infty)),\rho)$ respectively, (\ref{eq:MAMWIDAREJointCong}) is proved.

Secondly, we prove that (\ref{eq:MAMWIDAREAnnealCong}). By
(\ref{eq:PhiANPHiA}), for any open set $G\subset D_K[0,\infty)$, we
have
\begin{equation}
\liminf_{N\rightarrow \infty} \Phi_N(\hat{A}_N(\hat{\omega}))(G)\ge
\Phi(\hat{A}(\hat{\omega}))(G), \mbox{ $\hat{Q}$-a.s. }
\end{equation}
 Then by Fatou Lemma, we get
\begin{equation}
 \liminf_{N\rightarrow
\infty} \int \Phi_N(\hat{A}_N(\hat{\omega}))(G)
\hat{Q}(d\hat{\omega})\ge \int
\Phi(\hat{A}(\hat{\omega}))(G)\hat{Q}(d\hat{\omega}),
\end{equation}
which implies that
\begin{equation}\liminf_{N\rightarrow \infty}
\int \Phi_N(A_N(\omega_N))(G) Q_N(d\omega_N)\ge \int
\Phi(A(\omega))(G)Q(d\omega).
\end{equation}
 Since $D_K[0,\infty)$
is separable, (\ref{eq:MAMWIDAREAnnealCong}) is proved.

\appendix

\section{Weak Convergence Criteria for Time Inhomogeneous Markov
Processes}\label{se:WCCFTInhomMP}

In this section we state the weak convergence criteria for time
inhomogeneous Markov processes. These criteria are concerned with
martingale problems with time-dependent generators. We introduce
some notations from Either and Kurtz (1986). For $n=1$, 2, $\cdots$,
let
  $\{\mathscr{G}^n_t\}$ be a complete filtration, and let $\mathscr{L}_n$ be the space of
  real-valued $\{\mathscr{G}^n_t\}$-progressive
   processes $\xi$ satisfying
$$\sup_{0\le t\le T}E[|\xi(t)|]<\infty$$
for each $T>0$. Let $\hat{\mathscr{A}}_n$ be the collection of pairs
$(\xi, \varphi)\in \mathscr{L}_n\times \mathscr{L}_n$ such that
$$ \xi(t)-\int_0^t\varphi(s)ds$$
is a $\{\mathscr{G}^n_t\}$-martingale.

{\prop\label{MarThmInhom} Let $(E,\tilde{r})$ be a Polish space. Let
$\{G(s), 0\le s<\infty\}$ be a family of operators
 on $\bar{C}(E)$. Suppose that there exists a countable set $\Gamma_1\subset [0, \infty)$, such that for each
  $s\notin \Gamma_1$, $\{G(s), 0\le s<\infty\}$ has a common domain denoted by $\mathscr{D}(G)$ and for each
   $f\in \mathscr{D}(G)$, $G(s)f\in \bar{C}(E)$ for $s\notin \Gamma_1$, and $\|G(s)f\|$ is bounded for $s\in [0, T]\setminus\Gamma_1 $
   for any $T>0$. Suppose that there is an algebra $C_a$ contained
in the closure of $\mathscr{D}(G)$
     (in the sup norm) which separates points. Suppose that the $D_E[0,\infty)$ martingale problem for $(G, \nu)$
      has at most one solution, where $\nu\in \mathscr{P}(E)$. Suppose for each $n\ge 1$, $X_n$ is a
       $\{\mathscr{G}^n_t\}$-adapted process with sample paths in $D_E[0,\infty)$. Suppose $P (X_n(0))^{-1}\Rightarrow \nu$
        and the compact containment condition holds. Suppose $M\subset \bar{C}(E)$ is separating. Then
        condition (\textbf{a}) implies that there exists a solution $X$ of the $D_E[0,\infty)$ martingale problem for $(G, \nu)$, and
          $X_n\Rightarrow X$:

(\textbf{a}) There exists a countable set $\Gamma_2\subset [0,
\infty)$ such that for each $f\in \mathscr{D}(G)$,
 and $T>0$, there exists $(\xi_n, \varphi_n)\in \hat{\mathscr{A}}_n$, such that
\begin{align}
\label{eq:xinbd} &\sup_n \sup_{0\le s\le T,s\notin \Gamma_1 }E[|\xi_n(s)|]<\infty,   \\
\label{eq:varphinbd} &\sup_n \sup_{0\le s\le T,s\notin \Gamma_1}E[|\varphi_n(s)|]<\infty,   \\
\label{eq:xincon}\lim_{n\to \infty}E[&(\xi_n(t)-f(X_n(t)))\prod_{i=1}^{k}h_i(X_n(t_i))]=0, \\
\label{eq:varphincon}\lim_{n\to
\infty}E[&(\varphi_n(t)-G(t)f(X_n(t)))\prod_{i=1}^{k}h_i(X_n(t_i))]=0,
\end{align}
for all $k\ge 0$, $0\le t_1<t_2<\cdots<t_k\le t\le T$ with $t_i$,
$t\notin \Gamma_1\cup\Gamma_2$, and $h_1$, $\cdots$, $h_k\in M$, and
\begin{align}
\label{eq:xinfcon0} &\lim_{n\to \infty}E[\sup_{t\in \tilde{Q}\cap [0, T]}|\xi_n(t)-f(X_n(t))|]=0,   \\
\label{eq:varphinbdlpnorm} &\sup_n E[\|\varphi_n\|_{p,T}]<\infty,
\mbox{ for some $p\in (1,\infty]$},
\end{align}
where $\tilde{Q}$ is a countable and dense subset of $R$,
$\|h\|_{p,T}=[\int_0^T|h(t)|^p dt]^{1/p}$ if $p<\infty$;
$\|h\|_{\infty, T}=\mbox{ess}\sup_{0\le t\le T}|h(t)|$. }

The conditions of the following two corollaries are more convenient
to be verified for our multiagent models.

{\cor\label{InhomMarkovGntGtCon} Suppose in Proposition
\ref{MarThmInhom} that for each $n\ge 1$, $X_n$ is a
$\{\mathscr{G}^n_t\}$-adapted process with sample paths in
$D_E[0,\infty)$ and generator $\{G_n(s), 0\le s<\infty\}$ on
$\bar{C}(E)$. Suppose also for each $n\ge 1$ there exists a
countable set $\Gamma^n\subset [0, \infty)$, such that for each
$s\notin \Gamma^n$, $\{G_n(s), 0\le s<\infty\}$ has a common domain
denoted by $\mathscr{D}(G_n)$. Then condition (\textbf{b}) implies
that there exists a solution $X$ of the $D_E[0,\infty)$ martingale
problem for $(G, \nu)$, and $X_n\Rightarrow X$:

(\textbf{b}) For each $n\ge 1$, $\mathscr{D}(G_n)=\mathscr{D}(G)$
and there exists a countable set $\Gamma_2\subset [0, \infty)$
 such that for each $f\in \mathscr{D}(G)$, and $T>0$,
\begin{equation*}
\sup_n\sup_{0\le s\le T, s\notin
\Gamma_1\cup\Gamma^n}\|G_n(s)f\|<\infty,
\end{equation*} where
$\|\cdot\|$ is the sup norm on $B(E)$, and
\begin{equation*}
\lim_{n\to \infty}\sup_{q\in E}|G_n(t)f(q)-G(t)f(q)|=0,
\end{equation*} for any $t$ satisfying $t\notin
\Gamma_1\cup\Gamma_2\cup_{n=1}^{\infty}\Gamma^n$ and $0\le t\le T$.
}

{\cor\label{InhomMarkovChainCon} Suppose in Proposition
\ref{MarThmInhom} that $X_n=\eta_n(Y_n([\alpha_n\cdot]))$
 and $\{\mathscr{G}^n_t\}=\{\mathscr{F}^{Y_n}_{[\alpha_n t]}\}$, where $\{Y_n(k), k=0,1,2,\cdots \}$ is a time inhomogeneous Markov chain in a metric space $E_n$ with transition functions $\mu_{n,k}(x,\Gamma)$, $\eta_n: E_n\mapsto E$ is Borel measurable, and $\alpha_n\rightarrow \infty$ as $n\rightarrow \infty$. Define $T_{n,k}: B(E_n)\mapsto B(E)$ by
\begin{equation*}
T_{n,k}f(x)=\int f(y)\mu_{n,k}(x,dy),
\end{equation*}
 and let
$G_{n,k}=\alpha_n(T_{n,k}-I)$. Then condition (\textbf{c}) implies
that there exists a solution $X$ of
 the $D_E[0,\infty)$ martingale problem for $(G, \nu)$, and $X_n\Rightarrow X$:

(\textbf{c}) There exists a countable set $\Gamma_2\subset [0,
\infty)$ such that for each $f\in \mathscr{D}(G)$, and
 $T>0$,
\begin{equation}\label{eq:essANbd}
\sup_n\sup_{0\le t\le T, t\notin \Gamma_1\cup\Gamma_2}\sup_{q\in
E_n}|G_{n,[\alpha_n t]}f\circ \eta_n(q)|<\infty,
\end{equation} and
 \begin{equation}\label{eq:GnGtEnNormcon}
 \lim_{n\to \infty}\sup_{q\in E_n}|G_{n,[\alpha_n t]}f\circ \eta_n(q)-G(t)f\circ \eta_n(q)|=0,
 \end{equation}
for any $t$ satisfying $t\notin \Gamma_1\cup\Gamma_2$ and $0\le t\le
T$. }

{\rem\label{rem:GenAndSimp} \textbf{1)} Assume that $D\subset
\mathscr{D}(G)$, and $C_a\subset D$
 where $C_a$ is a subalgebra of $\bar{C}(E)$ which separates points.
 If we can prove that the uniqueness of $D_E[0,\infty)$ martingale problem for $(G,
 \nu)$ holds for functions in $D$ instead of the common domain $\mathscr{D}(G)$, then we can replace
  $\mathscr{D}(G)$ by $D$ in the conditions (\textbf{a}), (\textbf{b}) and (\textbf{c}). \textbf{2)} We can simplify the condition (\textbf{c}) of Corollary
\ref{InhomMarkovChainCon} into the following version: there exists a
countable set $\Gamma_2\subset [0, \infty)$ such that for each $f\in
\mathscr{D}(G)$, and  $T>0$,
\begin{equation}\label{eq:GnGtUniformEnNormcon}
\lim_{n\to \infty}\sup_{0\le t\le T, t\notin
\Gamma_1\cup\Gamma_2}\sup_{q\in E_n}|G_{n,[\alpha_n t]}f\circ
\eta_n(q)-G(t)f\circ \eta_n(q)|=0.
\end{equation}
This condition is stronger than the condition (\textbf{c}), since
(\ref{eq:GnGtUniformEnNormcon}) implies (\ref{eq:essANbd}) and
(\ref{eq:GnGtEnNormcon}). }

\section*{Acknowledgements}

The author would like to thank Professor Donald A. Dawson for his
valuable suggestions of the formulation of the multiagent models.

\end{document}